\documentclass[11pt]{amsart}
\usepackage{amsmath,amssymb,enumerate}
\usepackage{epsf,epsfig,amsfonts,graphicx,color,url}
\usepackage{pgf,tikz}
\usepackage{ulem}

\usepackage[all]{xy}
\usepackage{filecontents}

\def\R{\mathbb{R}}

\def\D{\mathcal{D}}
\def\J{\mathcal{J}}

 \newcommand{\dd}[2]
{
{{\partial #1}   \over {\partial #2}}
}

\numberwithin{equation}{section}

\newtheorem{definition}{Definition}[section]

\newtheorem{theorem}[definition]{Theorem}
\newtheorem{proposition}[definition]{Proposition}
\newtheorem{corollary}[definition]{Corollary}
\newtheorem{remarkth}[definition]{Remark}

\renewcommand{\emph}[1]{{\bfseries\itshape{#1}}}

\numberwithin{figure}{section}

\def\sR{subRiemannian } 
\def\on{orthonormal }

\begin{document}
\title[Higher elastica.   ]{ Higher Elastica:  geodesics in the Jet Space}  
\author[A.\ Bravo-Doddoli]{Alejandro\ Bravo-Doddoli}
\address{Alejandro Bravo: Dept. of Mathematics, UCSC,
1156 High Street, Santa Cruz, CA 95064}
\email{Abravodo@ucsc.edu}
\begin{abstract}
Carnot groups are  \sR manifolds.  As  such they admit geodesic flows,
which are left-invariant Hamiltonian flows on their cotangent bundles.  Some of these   flows  are integrable.  Some are not.
The space of k-jets for real-valued functions on the real line forms a Carnot group of dimension $k+2$. 
We show that its geodesic flow is  integrable and that its geodesics 
generalize     Euler's elastica, with the case  $k=2$ corresponding to the   elastica, as shown in \cite{ardentov2013conjugate}.
\end{abstract}
\thanks{}

\maketitle

\section{Introduction}

The space of k-jets of real functions of a single real variable, denoted here by $\J^k$,  is a $k+2$-dimensional manifold endowed with a canonical rank 2 distribution,
by which we mean   a linear sub-bundle of its tangent bundle.   This distribution is framed
by two vector fields, denoted $X_1, X_2$ below, whose iterated Lie brackets give
$\J^k$ the structure of a Carnot group.  Declaring $X_1$ and $X_2$ to be \on
  endows $\J^k$ with the structure of a \sR manifold,  one which is (left-) invariant under the  Carnot group multiplication.
Like any \sR structure, the    cotangent bundle $T^* \J^k$ is   endowed with  a Hamiltonian system whose underlying Hamiltionian $H$ 
is that whose solution  curves project to the  \sR geodesics on $\J^k$.   We call this  Hamiltonian system
the geodesic flow on $\J^k$.  
\begin{theorem}
\label{thm1} 
The geodesic flow for the \sR structure on $\J^k$ is integrable.
\end{theorem}

  $\J^1$ is isometric to the  
Heisenberg group where  this theorem is well-known see \cite{hakavuori2018blowups} and \cite{montgomery2002tour}.  $\J^2$ is isometric to  the Engel group and
Ardentov and Sachkov showed that its 
 \sR geodesics    correspond to Euler elastica.  Their  result  inspired our next theorem. 

\begin{figure}[h]
\centering
{\includegraphics[width=3cm]{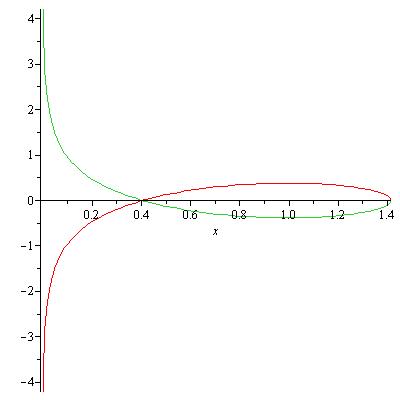}} 
\qquad
{\includegraphics[width=3cm]{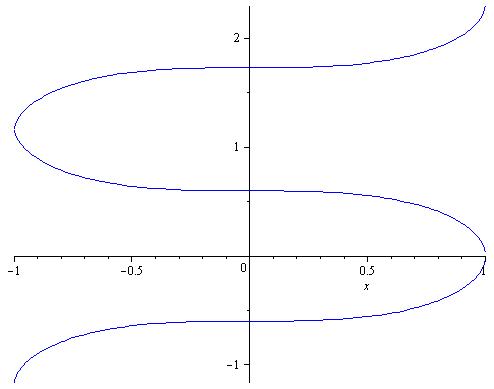}}
\qquad
{\includegraphics[width=3cm]{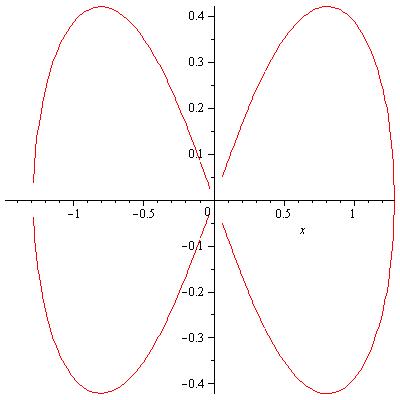}}
\caption{Some clasic solutions of the Elastica equation generated by $F_2(x) = \frac{x^2}{a^2}-\alpha$, on the left the convict curve $\alpha=1$, in the center the pseudo-sinusoid $\alpha =0$ and on the left the pseudo-lemniscate with $\alpha =.65222...$.} \label{fig:Andrei-solution }
\end{figure}

$\J^k$   comes with a projection $\Pi: \J^k \to \R^2 = \R^2 _{x, u_k}$
onto the Euclidean plane  which projects  the frame $X_1, X_2$
projects onto the standard coordinate frame
$\dd{}{x}, \dd{}{u_k}$ of $\R^2$.   (See SETUP below for the meaning of the coordinates.)   As a consequence, 
a horizontal curve   $\gamma$ in $\J^k$ is parameteried by
(\sR) arclength if and only if  its planar projection $\Pi \circ \gamma$ to $\R^2$ is  parameterized by arclength.
We will characterize geodesics on $\J^k$ in terms of
their planar projections.  As alluded to already, Ardentov and Sachkov, \cite{ardentov2013conjugate},  proved that when $k=2$ the planar projections of  geodesics
are Euler elastica.  These elastica have ```directrix'' the $u_2$-line, the line orthogonal to the $x$-axis. 
 There are many ways
to characterize Euler's  elastica, see \cite{bryant1986reduction}, \cite{bonnard2001stratification} and \cite{lawden2013elliptic}, \cite{jurdjevic1997geometric}.  The one  we will  use  is as follows.
Take a planar curve  $c(s) = (x(s), y(s))$  and consider its  curvature $\kappa = \kappa(s)$,
where $s$ is arclength.  Then the  curve $c$  is an Euler elastica with directrix a   line parallel to the y-axis
if and only if   $\kappa(s)) = P(x(s))$ for $P(x)$ some {\bf linear}  polynomial  in $x$ -- that is
$P(x) = a x + b$ for some constants $a$ and $b$. 
See FIGURE \ref{fig:Andrei-solution }. 

\begin{theorem} 
\label{thm2} Let $\gamma : I \to \J^k$ be a \sR geodesic parameterized by arclength $s$ and
$\pi \circ \gamma = c(s) = (x(s), u_k (s))$   its planar projection.
Let $\kappa$ be the curvature of $c$.  Then  $\kappa(s) = p(x(s))$
for some degree $k-1$-polynomial $p(x)$ in the coordinate $x$.   Conversely,
  any plane curve $c(s)$ in the $(x, u_k)$ plane which is  parameterized by arclength $s$  and   whose curvature $\kappa(s)$
equals $p(x(s))$ for some  polynomial $p(x)$ of degree at most  $k-1$   in $x$
is the projection of such  a \sR geodesic.
\end{theorem}

{\bf Example} For  the case  $k=1$ of the Heisenberg group the theorem  asserts that
$\kappa = P(x)$ where  $P$ a degree $0$ polynomial -- i.e. a constant function.
The only   curves having constant curvature  are lines and circles, and these  are well-known to be the
projections of the Heisenberg geodesics.

\section{set-up}

The  k-jet  of a smooth  function $f: \R \to \R$ at a point $x_0 \in \R$ is its kth order Taylor expansion
at $x_0$.  We will 
encode this  k-jet as a   $k+2$-tuple of real numbers as follows: 
\begin{equation} (j^k f )(x_0)  = (x_0, f^k(x_0),f^{k-1}(x_0),\dots,f'(x),f(x_0))   \in \R^{k+2}
\label{k-jet}
\end{equation}
As $f$ varies over  smooth functions and $x_0$ varies over $\R$,  these  $k$-jets
sweep out the   $k$-jet space, denoted   by $\mathcal{J}^k$.  $\J^k$   is diffeomorphic to $\R^{k+2}$
and   its   points are coordinatized according to
 $$(x,u_k,u_{k-1},\dots,u_1,y) \in \mathbb{R}^{k+2}: = \mathcal{J}^k.$$
Recall that if $y = f(x)$ then $u_1= dy/dx$ while  $u_{j+1} = du_j /dx$, $j \ge 1$.  Rearranging these equations into
$dy = u_1 dx, du_j = u_{j+1} dx$ we see that   $\J^k$  is endowed with a natural rank 2 distribution $\mathcal{D}\subset T\mathcal{J}^k$
characterized by the $k$ Pfaffian equations
 \begin{eqnarray*}
 u_1dx-dy & = & 0 \\
 u_2dx-du_1 & = & 0 \\
   \vdots & = & \vdots \\
 u_kdx-du_{k-1} & = & 0 
  \end{eqnarray*}

 The typical integral curves of $\D$ are the k-jet curves $x \mapsto (j^k f) (x)$. 
 In addition to these integral curves we have a distinguished family of 
 curves which arise by varying only the highest derivative $u_k$, and which are
 the integral curves of the vector field $X_2$ below (eq (\ref{frame}).  These latter curves are $C^1$-rigid
 in the sense of Bryant-Hsu,\cite{Bryant1993}, and they exhaust the supply of $C^1$-rigid curves.  
 
A \sR structure on a manifold consists of a non-integrable distribution together with a smoothly varying
family of inner products on the distribution.  We  have our distribution $\D$ on $\J^k$.  We arrive at our \sR structure 
by observing that
 $\D$ is  globally framed by the two vector fields 
 \begin{equation}
 \label{frame}
X_1 = \frac{\partial}{\partial x} + u_1 \frac{\partial}{\partial y} + \sum_{i=2}^{k} u_i \frac{\partial}{\partial u_{i-1}} \;\; \text{and} \;\; X_2 = \frac{\partial}{\partial u_k} 
\end{equation}
and then declaring   these two vector fields to be  \on.   
Now the restrictions of the one-forms  $dx, du_k$ to $\D$  form a global co-frame for $\D^*$ which is dual to our frame (eq (\ref{frame}). 
It follows that 
an equivalent way to describe our \sR structure is to say that its  metric is  $dx^2 + du_k ^2$ {\bf restricted to } $\D$.

 For the purposes of theorem \ref{thm2} the following alternative characterization of the \sR metric is crucial.
 Consider the projection
 $$\Pi:  \mathcal{J}^k \to \R^2 _{x, u_k };  \Pi(x,u_k,u_{k-1},\dots,u_1,y) = (x, u_k).$$
 Its fibers  are transverse to $\D$ and we have that $\Pi_* X_1 = \dd{}{x}, 
 \Pi_*X_2 = \dd{}{u_k}$, so our frame pushes down to the standard frame for $\R^2$  The metric on each two-plane $\D_p$,
 $ p \in \J^k$ is 
  characterized 
 by the condition that $d \Pi _p$ (which is just $\Pi$ since
 $\Pi$ is linear), restricted to  $\D_p$  is a linear isometry onto $\R^2$, where $\R^2$ is endowed with the standard metric
 $dx^2 + du_k^2$.  
 It follows immediately that the length of any horizontal path equals the length of its planar projection,
that $\Pi$ is a ``submetry'': $\Pi(B_r (p)) = B_r (\Pi (p))$, where $B_r (p)$ denotes the metric ball of radius $r$ about $q$,
and that the horizontal lift a  Euclidean line
in $\R^2$ is a geodesic in $\J^k$.

\subsection{Hamiltonian }  

Let $P_1, P_2:  T^* \J^k \to \R$  be the `power functions' of the vector fields $X_1,X_2$ above.
(REF: \cite{montgomery2002tour},  8 pg.) . In terms of traditional cotangent coordinates 
$(x, u_k, u_{k-1}, \ldots, u_1, y, p_x, p_k, p_{k-1}, \ldots, p_1, p_y)$ for $T^* \J^k$, with $p_i$ short for $p_{u_i}$
 we have
$$P_1 = p_x + u_1 p_y + u_2 p_1 + \ldots + u_k p_{k-1}; \qquad P_2 = p_k.$$
Then the Hamiltonian governing the \sR  geodesic flow on $\J^k$ is 
\begin{equation} H = \frac{1}{2} (P_1 ^2 + P_2 ^2 ).
\label{K.E.}
\end{equation}
See \cite{montgomery2002tour}; 8 pg.
If we want our geodesics to be parameterized by arclength then we set $H=1/2$,
and this we will do in what follows. 

{\bf Remark.}  [$C^1$-rigidity.]  The $u_k$ curves are $C^1$-rigid for $\D$, and
form what Liu-Sussmann christened as the  ``regular-singular'' curves for $\D$.
As such, they are geodesics for {\bf any} \sR metric $E dx^2 + 2F dx du_k + G du_k ^2$, restricted to $\D$. Such that $ds^2$ is positive definite, for $E, F, G$ any functions of the jet coordinates $(x,u_k, u_{k-1}, \ldots, y)$, regardless of whether or not they satsify the corresponding
(normal) geodesic equations. For our   metric
  each  $u_k$-curve is   indeed the projection to $\J^k$ of a solution to our $H$, so we do not go to
  extra effort to account for these abnormal geodesics.  (REF \cite{montgomery2002tour},  chapter 3).

 \section{Carnot Group structure} 

Under iterated bracket our  frame $\{ X_1, X_2 \}$ generates a $k+2$-dimensional nilpotent Lie algebra
which can be identified pointwise with the tangent space to $\J^k$.
Specifically, if
we write
$$X_3 = [X_2, X_1],  X_{4} =[X_{3},  X_1], \ldots , X_{k+2} =[X_{k+1},  X_1], 0 = [X_{k+2},  X_1],  $$
then we compute that
$$X_{k+2} = \dd{}{y},  X_{k+1} = \dd{}{u_1},  X_k= \dd{}{u_2}, \ldots,  X_3 = \dd{}{u_{k-1}}$$

and that all other Lie brackets $[X_i, X_j]$ are zero.  
The span of  the $X_i$ thus form a $k+2$-dimensional graded nilpotent Lie algebra
$$\mathfrak{g}_k = V_1 \oplus V_2 \oplus \ldots V_{k+1}, \;\;\; V_1  = span\{ X_1, X_2\} ,   V_i =span\{ X_{i+1}\} ,  1 < i \leq k+1.$$
Like any graded nilpotent Lie algebra, this algebra has an associated Lie group
which is a Carnot group $G$, and by using the flows of  the $X_i$ we can identify
$G$ with $\J^k$, and the $X^i$ with left-invariant vector fields on $G \cong \J^k$.

\section{ Integrability:  Proof of theorem \ref{thm1}}.

Our Hamiltonian $H$ is a left-invariant Hamiltonian on the cotangent bundle of a Lie group $G$.
We recall the general `Lie-Poisson'' structure for such Hamiltonian flows.    \cite{arnol2013mathematical} Appendix;
REF: local struc Poisson, \cite{marsden2013introduction} ch 4.

\begin{equation*}
\xymatrix{
& T^*G   \ar[dl]^{J_L}  \ar[dr] ^{J_R }&  \\
{\frak g^* _+}   & &  {\frak g^* _-} 
}\end{equation*}

The arrows $J_R, J_L$ are the momentum maps for the actions of $G$ on itself by  right and left translation, lifted to $T^*G$. 
The subscripts $\pm$ are for a plus or minus sign in front of the Lie-Poisson (=Kostant-Kirrilov-Souriau) bracket on $\frak g^*$.
$J_R$ corresponds to {\it left} translation back to the identity and realizes the quotient of $T^*G$ by the {\it left} action. 
$J_L$ corresponds to {\it right} translation of a covector back to the identity and forms the components of the momentum
map for {\it left} translation, lifted to the cotangent bundle.
 In our case,
${\frak g^* } = \R^{k+2}$ and 
$$J_R = (P_1, P_2, P_3, \ldots, P_{k+2})$$
with $P_i$ the power function associated to $X_i$, so that 
$$P_3 = p_{k-1}, P_4= p_{k-2}, \ldots , P_{k+2} = p_y$$
in terms of standard canonical coordinates as above.

When the Hamiltonian $H: T^*G \to \R$ is  left-invariant  it can be expressed as a function of the components
of $J_R$,  that is $H = h \circ J_R$ for some $h: {\frak g}^* \to \R$, and $H$  Poisson commutes
with every component of the   {\it left}  momentum map 
$J_L$, so that these left-components are   invariants.  $J_L$ and $J_R$ are related by
$J_L (g, p)  = (Ad_g ) ^* J_R (g,p)$ where we have written $p \in T_g ^* G$ and where $Ad_g$ is the
adjoint action of $g$.     

The reason underlying the integrability of our system is a   simple dimension count. 
\begin{proposition} If the generic  co-adjoint orbit of ${\frak g}^*$ is 2-dimensional
then the left-invariant Hamiltonian flow on $T^*G$ is integrable.
\end{proposition}

Recall that the symplectic reduced spaces for {\it left translation} action are the co-adjoint orbits,
for ${\frak g}^*_+$, and that $J_R$ realizes this symplectic reduction procedure, mapping each
$J_L^{-1}(\mu)$ onto the co-adjoint orbit through $\mu$.  The hypothesis of the Proposition asserts that the 
symplectic reduced spaces associated to the $G$-action  are zero or two  dimensional,
so, morally speaking, the system is automatically integrable by reasons of dimension count.     

{\sc Proof of Proposition.} We must produce  $n$ commuting
integrals in involution, where  $n = dim(G)$.  The hypothesis  asserts that  there are
$n-2$ Casimirs $C_1, \ldots , C_{n-2}$ for ${\frak g}^*$, these being the  functions whose 
common level sets  at a generic value define a generic  co-adjoint orbit.  These Casimirs
are a functional basis for the $Ad^*_G$ invariant polynomails on ${\frak g}^*$.
When viewed as functions on $T^*G$ via $C_i \circ J_R$
the Casimirs  Poisson commute with {\it any} left-invariant function on $T^*G$, and in 
particular with $H$ and with each other.   Thus, $H, C_1, C_2, \ldots, C_{n-2}$
yield $n-1$ integrals.  To get the last commuting integral take any component of 
$J_L$.  
QED.   

{\sc Proof of theorem  \ref{thm1}.} In order to use the proposition, we need to verify that the co-adjoint orbits
are generically 2-dimensional.      We have the Poisson brackets
\begin{equation}
\{ P_1, P_i \} = P_{i+1},  1 < i < k+1,\qquad \text{and} \qquad \{P_1,P_{k+2}\} = 0.
\end{equation}
with all other Poisson brackets $\{P_i, P_j \}, 1< i < j \le k+2$ being  zero.  Thus the Poisson tensor $B$  at a point 
$Z = (P_1, P_2, P_3, \ldots, P_{k+2}) \in {\frak g}^* _+$ is :
\begin{equation}
B:= \begin{pmatrix}
0 & Z_k  & 0\\
Z_k^t & 0 & 0  \\
0 & 0 & 0 \\
\end{pmatrix} \;\; \text{where} \;\; Z_k = ( P_3, \ldots, P_{k+2} ).
\end{equation}
which  has rank $2$ generically and rank $0$ if and only if $Z_k = 0$, i.e. if and only if  $P_i =0$ for $2<i$.
QED

Thanks to this information we now that the system has $k$ Casimir functions.  
\begin{theorem}\label{Casimir}
If $P_{k+2}\neq 0$ the function $C_i( P_2, P_3, \ldots, P_{k+2})$, with $i \in \{1,\dots,k\}$ which are given by $C_1(P_2, P_3, \ldots, P_{k+2}) = P_{k+2}$,  and for  $i>1$, we define as 
\begin{equation}\label{Casimir-f}
\begin{split}
C_i =&\; P_{k+2}^{i-1}P_{k+2-i} +\sum_{j=1}^{i-2}(-1)^{j} P_{k+2}^{i-(j+1)} P_{k+2-j} \frac{P_{k+1}^j}{(j)!} +(-1)^{i-1} \frac{P_{k+1}^i}{(i-2)!\;i}    ,      
\end{split}
\end{equation}
are constant of motion for the geodesic equations in $\mathcal{J}^k$, in others words they are Casimir functions\footnote{In the case $i=2$ the sum is empty}. 
\end{theorem}

\section{Integration and curvature: proof of theorem \ref{thm2}}

Hamilton's equations read $\frac{df}{dt} = \{ f, H\}$.  
With our Hamiltonian they expand out to $\frac{df}{dt} = \{f, P_1 \} P_1 + \{f,P_2 \}P_2$.
Returning to our coordinates $x, u_k$ we compute $\{ x, P_1 \} = 1,  \{x, P_2 \} = 0$,
$\{ u_k , P_1 \} = 0,  \{ u_k , P_2 \} = 1$ so that
\begin{equation}
\frac{ dx}{dt} = P_1 \qquad \text{and} \qquad \frac{du_k}{dt} = P_2.
 \label{R-H-E-0}
\end{equation} 
Thus $(P_1, P_2)$ are the components of the tangent vector to the plane curve $(x(s), u_k (s))$
obtained by projecting a geodesic to the plane.  If $H = 1/2$ then
this vector  is a unit vector, the parameter $t$ of the flow is arc-length $s$ and we can write 
\[ (P_1,  P_2 ,   \dot \theta ) = (\cos(\theta),\sin (\theta) ,\kappa) \]
Using $\{P_1, P_2 \} = P_3$ we see that the $P_1$,  evolve according to
\begin{equation} \label{R-H-E-1}
\begin{split}
\dot P_1 & = P_3 P_2 \\
\dot P_2 & = -P_3 P_1. \\
\end{split}
\end{equation}
But we also have $\dot P_1 = - \dot \theta \sin(\theta) = -\dot \theta P_2$ and
$\dot P_2 = + \dot \theta \cos(\theta) = + \dot \theta P_1$
from which we see that
$$-P_3 = \kappa.$$
Now for $k+2> i >2$ we have that $\{ P_i, P_1 \} = - P_{i+1},  \{ P_i, P_2 \} = 0$, and $\{ P_j, P_{k+2} \} = 0  $ for all $j$ so that 
\begin{equation} \label{R-H-E-2}
\begin{split}
\dot P_3 &= -P_1 P_4  ,\\
\dot P_4 &= -P_1 P_5  ,\\
\vdots \\
\dot P_{k+1} &= - P_1 P_{k+2} ,\\
\dot P_{k+2} &= 0.  \\
\end{split}
\end{equation}

{\sc Proof of theorem \ref{thm2}.} Consider a geodesic $\gamma$  and an  arc of the  geodeisc for which   $\dot x \ne 0$.
Instead of arclength $t = s$    use $x$ to parameterize this arc.
From eq (\ref{R-H-E-0}) we have, along this arc, that  
$$\frac{d}{dx} = \frac{1}{P_1} \frac{d}{ds}$$
so that the equations for the evolution of $P_3, P_4, \ldots, P_{k+2}$ along the curve become
$$\frac{dP_3}{dx} = - P_4$$
$$\frac{dP_4}{dx} = - P_5$$
$$\ldots $$
$$\frac{dP_{k+1}}{dx} = - P_{k+2}$$
$$\frac{dP_{k+2}}{dx} =0$$
These equations can  be summarized by
$$\frac{d ^k P_3}{dx ^k} = 0$$
which asserts that the curvature $P_3$ of the projected curve $c = \pi \circ \gamma$,
is  a  polynomial $p(x)$ of degree $k-1$ in $x$,    at least along our  arc.
Finally, since $\gamma$ is an analytic function of $s$, so are $c(s)$ and $\kappa(s)$,
so that if $\kappa(s)$ enjoys a relation
$\kappa(s) = p(x(s))$ along some subarc of $c(s)$, it enjoys this same  relation
everywhere along $c$.

To prove the converse,  first  consider  a  general smooth curve in the $x-u$ plane
along which $dx > 0$.  We can parameterize the curve either by arc-length
$s \mapsto (x(s), u(s))$ or as graph, $u = u(x)$.  
Define the function $F(x)$, with   $-1 \le F(x) \le 1$ by way of relating the two parameterizations:  
\begin{equation}
\label{F1}(\dot x , \dot u):= (\frac{dx}{ds}, \frac{du}{ds})  = (\sqrt{1- F(x)^2 },  F(x))
\end{equation}
so that  $dx = \sqrt{1 - F(x)^2} ds$ and 
$$u'(x):= \frac{du}{dx}  = \frac{F(x)}{\sqrt{ 1- F(x)^2}}$$
It follows    that 
$$u'' = \frac{F'}{(1- F(x)^2)^{3/2}}. $$
Now the curvature of our  curve, when viewed as a graph, is well known to be  
$$\kappa(x) = \frac{u'' (x)}{(1 + u'(x)^2)^{3/2}}$$
and we have 
$$(1 + u'(x)^2) = \frac{1}{1 - F(x)^2}$$
from which we conclude that
\begin{equation}
\label{Fkappa}
\kappa = F'(x).
\end{equation}

To finish the proof,  suppose that  we are given a curve $c$ in the   $x-u$ plane, with $u = u_k$,   whose  curvature $\kappa$  is a degree $k-1$ polynomial in $x$.
Define  $F(x)$ by  eq (\ref{F1})  along an arc of $c$ for which $dx > 0$
From eq (\ref{Fkappa})  we know that  $F$ is an anti-derivative of $\kappa$ and so a polynomial of degree $k$ in $x$.
(The constant term in the integration $F(x) = \int ^x \kappa   dx$ is fixed by
choosing any point  $(x_*, u_*) = (x(s_*), u(s_*)) $ along  $c$  for which $dx/ds > 0$ so that $-1 < du/ds|_{s = s_*} < 1$
and setting  $F(x_*) =du/ds|_{s = s_*}$.)
By the preceding analysis,  $c$   has curvature $\kappa(x(s))$ along the entire arc $dx > 0$  of our  curve which contains $(x_*, u_*)$.
Moreover $u' (x) =  F(x)/\sqrt{1- F(x)^2}$.
Set 
\begin{equation}
\label{PtoF1}
( P_1 (x) , P_2 (x)) := ( \sqrt{1-(F(x) )^2 },F(x)) ,, 
\end{equation}
\begin{equation}
\label{PtoF2} P_3(x) := F'(x),   \;\; \text{and}\;\; P_{i+2} := (-1)^{i}\frac{d^i F}{d^i x}(x), i >1.
\end{equation}
View the $P_i$ as    momentum functions.
Reparameterize the momentum functions  by  $s$  using $dx/ds = P_1 (x)$.  Then we  
 verify that the $P_i$   satisfy   \ref{R-H-E-1} and \ref{R-H-E-2},
so that the horizontal curve   $\gamma(x(s))$ in $\mathcal{J}^k$ with  these momenta  satisfies the geodesic equations and   projects on our given curve $c$.

QED

\begin{corollary}
Suppose that  the momentum functions $P_i$ are related to the degree $k$ polynomial 
$F(x)$ as per equations (\ref{PtoF1}, \ref{PtoF2}) and that $H = 1/2$.   Then a critical point  $x_0$  of  $F(x)$  corresponds
to a  relative equilibrium for the reduced equations  \ref{R-H-E-1} and \ref{R-H-E-2} if and only if $F(x_0) = \pm 1$.
\end{corollary}

{\sc Proof.} The equilibria of equations (\ref{R-H-E-1}) and (\ref{R-H-E-2}) are the points with  $P_1 = 0$ and $P_{3} = 0$, as long as $H \ne 0$.  If   $H=\frac{1}{2}$, the condition
$P_1 = 0$ forces   $P_2 = \pm 1$  but $P_2 = F(x)$.  Finally  $P_3 = F'(x)$.

\section{Structure of higher Elastica}

As we see in the last prove we have an option to select a primitive of $p(x)$, then given $F_k(x) =  \int p(x)dx$ the dynamics is trivial when $F_k^{-1}([-1,1])$ is empty or isolated points, i.e. $F_(x(s)) $ is constant for all $s$. Then we can take $F_k^{-1}([-1,1])$  as follow 
\begin{equation*}
F_k^{-1}([-1,1]) := \cup_{I} [x_0^i,x_1^i] \qquad \text{where}  \qquad F_k(x^i_0) = \pm 1 , \quad F_k(x^i_1) = \pm 1,
\end{equation*} 
$x_0^i<x_1^{i}\leq x_0^{i+1}<x_1^{i+1}$
and the condition that $F_k(x)\neq \pm 1$ if $x \in (x_0^i,x_1^i)$, note all the possibilities. Choose  $[x_0^i,x_1^i] =[x_0,x_1]$, by the last corollary the points $x_0$ and $x_1$ are equilibrium points if and only if they are critical points of the function  $F_k(x)$, then it takes infinite time to arrive to them.

\begin{theorem}
The curve $(x,u)$ with curvature $k((x(s)) = p(x)$ is bounded in the $x$-direction, and generically the curve is periodic in $x$ and the period   $L$ given by
\begin{equation*}
L := \int_{x_0}^{x_1} \frac{2dx}{\sqrt{1-F_k^2(x)}} \;\; \text{we also define }\;\; \tau = \int_{x_0}^{x_1} \frac{ 2F_k(x)dx}{\sqrt{1-F_k^2(x)}}.
\end{equation*}
Finally, we have that $u(s+L) = u(s) + \tau$.
\end{theorem} 
Let $x_0$ be a regular point, we will answer the question how to extend the curve $c(s)$ as a function of $x$ such that its lift is a smooth solution of geodesics equation,  set   $(P_1,P_2)=(\cos\theta, \sin\theta)$ and $\dot{\theta} = p(x)$ since $F_k(x_0) = \pm1$ define $\theta(x_0) = \pm \frac{\pi}{2}$ and  $\dot{\theta}(x_0)\neq 0$, then $P_1$ has a change of sign, while, $P_2$ does not change. Therefore if $x(s_0) = x_0$ we define
\begin{equation}
(\dot{x},\dot{u}) = \begin{cases}   (\pm\sqrt{1-F_k^2(x)},F_k(x))\qquad \text{if} \qquad s_0-\frac{L}{2} \leq s \leq s_0, \\
                                     (\mp \sqrt{1-F_k^2(x)},F_k(x))\qquad \text{if} \qquad  s_0 \leq s \leq s_0 + \frac{L}{2}. 
                      \end{cases}
\end{equation}
Therefore, the curve stays in the interval $[x_0,x_1]$, same with $x_1$.  If both are regular point, we can read the equation $P_1 = \sqrt{1-F_k(x)}$ as the restriction
$P_1(x)|_{\{H=\frac{1}{2},C_1,\dots,C_k\}}$ and consider action function $I$ given by the area under the graph $\sqrt{1-F_k(x)}$ going from $x_0$ to $x_1$ and the area of $-\sqrt{1-F_k(x)}$ going from $x_1$ to $x_0$, i.e.
\begin{equation*}
\mathcal{I}(H=\frac{1}{2},C_1,\dots,C_k) := 2\int_{x_0}^{x_1} \sqrt{1-F_k(x)} dx.
\end{equation*}
Finally, the period is given by $\frac{\partial \mathcal{I}}{\partial H}|_{\{H=\frac{1}{2},C_1,\dots,C_k\}} = L$, (see \cite{arnol2013mathematical} chapter 10). The period goes to infinite when $x_0$ or $x_1$ are critical points, much like the very well known homolinic connection of a pendulum.

Let us consider $(x_0,u_*)$ the initial point of the curve and  $x(s) \in [x_0,x_1]$ and  $2s \leq L$, then 
\begin{equation*}
\begin{split}
u(s) +  \tau &= \int_{x_0}^{x(s)} \frac{F_k(x)dx}{\sqrt{1-F_k^2(x)}} + \int_{x_0}^{x_1} \frac{2F_k(x)dx}{\sqrt{1-F_k^2(x)}} + u_*, \\
                        & = (\int_{x_0}^{x(s)}  +\int_{x(s)}^{x_1} + \int_{x_0}^{x_1}+\int_{x_0}^{x(s+L)} ) \frac{F_k(x)dx}{\sqrt{1-F_k^2(x)}}  
+ u_* = u(s+L).\\
\end{split}
\end{equation*}
Where again we use the fact that $\int_{x_0}^{x_1}\frac{F_k(x)dx}{\sqrt{1-F_k^2(x)}} = \int_{x_1}^{x_0}\frac{F_k(x)dx}{-\sqrt{1-F_k^2(x)}}$. QED

Here, we have three cases;
\begin{itemize}
\item   Periodic case - 
$p(x_0) \neq 0 $ and $p(x_1) \neq 0 $
\item Asymptotic behavior to one line -
 $p(x_0) = 0 $ and $p(x_1) \neq 0 $ or $p(x_0) \neq 0 $ and $p(x_1) = 0 $   .
\item Asymptotic behavior to two line -
$p(x_0) = 0 $ and $p(x_1) = 0 $ .
\end{itemize}

\subsection{General Convict curve}
\begin{figure}\label{fig:geneal-convict}
    \centering
    {{\includegraphics[width=3cm]{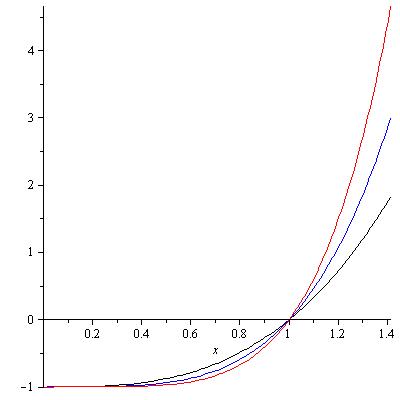} }}%
    \qquad
    {{\includegraphics[width=3cm]{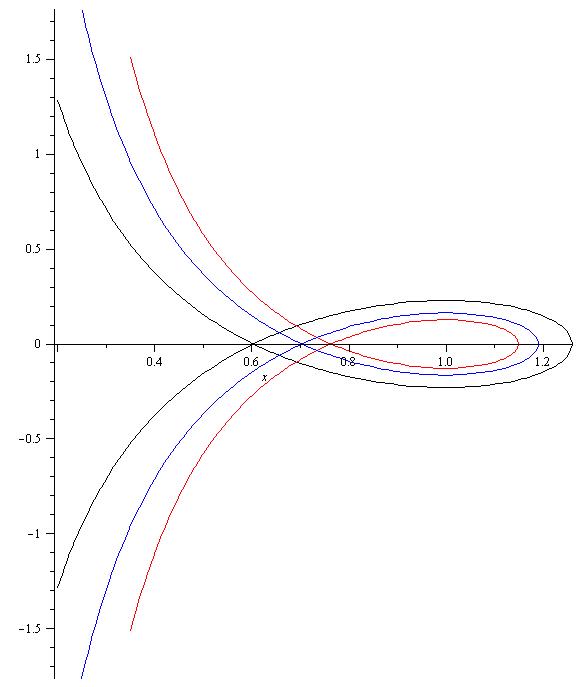} }}%
    \caption{On the left side we see  $F_k(x)$, while, on the right side we have the curves in the plane $(x,u)$ for $k=3,4,5$ .}%
\end{figure}

The elastica equation has a distinguished solution which we call the  \textbf{Euler Kink}. Other names for it are the Euler soliton or  Convict's curve.  See figure  and see \ref{fig:Andrei-solution }, see \cite{ardentov2013conjugate}, \cite{lawden2013elliptic}, \cite{jurdjevic1997geometric} and \cite{bryant1986reduction}.
We define the \textbf{Convict's curve} at the level $k$ in the sense that the curvature of the curve $(x,u_k(x))$ it is always proportional to $x^{k-1}$. See figure \ref{fig:geneal-convict} 
\begin{theorem}
If $1<k$ then the level $k$ has a convict curve.
\end{theorem}
Consider the polynomial $F_{k}(x) = \frac{x^{k}}{a^k} -\alpha$.  
Set $x =a\sqrt[k]{\alpha+\cos(t)}$ , we can find the next expressions
\begin{equation*}
\begin{split}
 u(t)  = \int_{t_0}^{t_1} \frac{ \cos(t) dt }{(\alpha+\cos(t))^{\frac{k-1}{k}}},&  \;\;\; t(x)  = \frac{a}{k} \int_{t_0}^{t_1} \frac{  dt }{(\alpha+\cos(t))^{\frac{k-1}{k}}}.  \\
\end{split}
\end{equation*}
The case $k=2$ is the classic solution for Elastica equation, see \cite{jurdjevic1997geometric} page 436. If $\alpha =1$, then we have the explicit expression
\begin{equation*}
\begin{split}
u(x) & = \int^x_{\sqrt[k]{2}} \frac{x^{\frac{k}{2}}dx}{\sqrt{2-x^{k}}} - \frac{2}{k\sqrt{2}} \ln(\frac{\sqrt{2}-\sqrt{2-x^k}}{x^\frac{k}{2}}), \\
 t(x) &  = \frac{1}{\frac{3}{2}\sqrt{2}} \ln(\frac{\sqrt{2}-\sqrt{2-x^k}}{x^\frac{k}{2}}) .\\
\end{split}
\end{equation*}
We can find a explicit second order ODE for $\dot{\theta}$,
\begin{equation*}
\begin{split}
p_{u_{k-1}} = \frac{\partial P}{\partial x}(x) = kx^{k-1} \;\;\; \text{and} \;\;\; \frac{\dot{\theta}^2}{k^2a^{2k}} = (\cos\theta + \alpha)^{\frac{2(k-1)}{k}}.\\
\end{split}
\end{equation*} 
In the case $k=2$ is the pendulum equation define in \cite{ardentov2013conjugate}, the ODE equation can be extend to $k=1$. Also, $\alpha = 1$ degenerate a system with a degenerated equilibrium point with a homoclinic connection.

\begin{figure}\label{fig:odd-graph}
    \centering
    {{\includegraphics[width=3cm]{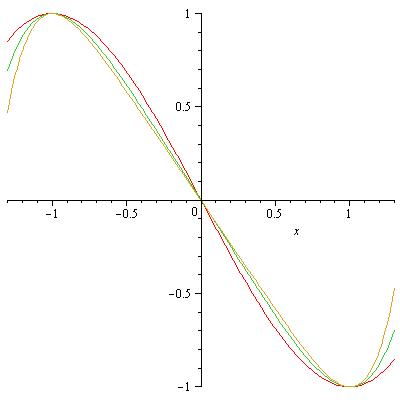} }}%
    \qquad
    {{\includegraphics[width=3cm]{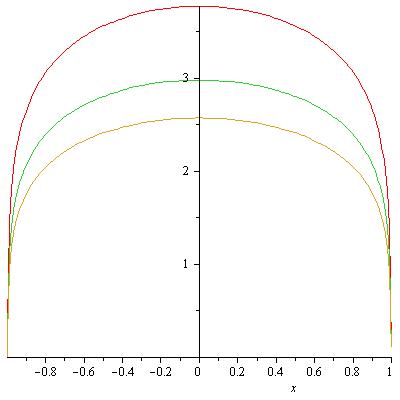} }}%
    \caption{On the left side we  see  $F_k(x)$, while, on the right side we have the curves in the plane $(x,u)$  for $k=3,5,7$.}%
\end{figure}

\subsection{Infinite geodesic graph.}
Here, we will define a infinite geodesic graph like the curve whose projection to the plane $(x,u)$ is always a graph of $x$. See figure \ref{fig:odd-graph} and \ref{fig:even-graph}.

\begin{theorem}
If $k>2$ then $\mathcal{J}^k$ has a geodesic graph.
\end{theorem}

We split in the even and odd case:
\begin{figure}\label{fig:even-graph}
    \centering
    {{\includegraphics[width=3cm]{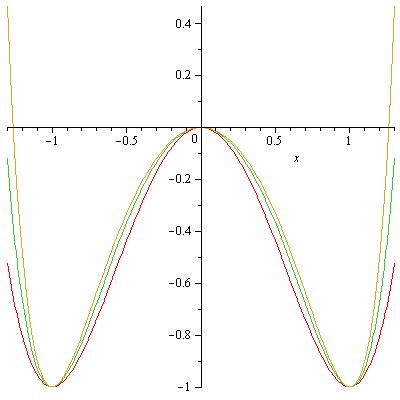} }}%
    \qquad
    {{\includegraphics[width=3cm]{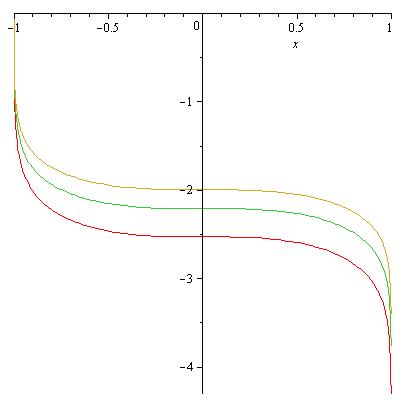} }}%
    \caption{On the left side we see  $F_k(x)$, while, on the right side we have the curves in the plane $(x,u)$ for $k=4,6,8$.}%
\end{figure}
\begin{itemize}
\item Consider the polynomial $F_{2k+1}(x) = -\frac{ x^{2k+1} -(2k+1)x }{2k}$ with $x \in [-1,1]$. Since the point $x = \pm 1$ give us $F_{2k+1}(\pm 1) =  \pm 1$ and $\frac{\partial F}{\partial x}(\pm 1) = 0$  they are equilibrium points.
\item Consider the polynomial  $F_{2k}(x) = -\frac{x^{2k}-k x^2}{1-k} $, with $x \in [-1,1]$. Since the point $x = \pm 1$ give us  $F_{2k}(\pm1) = -1$, and $\frac{\partial F}{\partial x}(\pm 1) = 0$. 
\end{itemize}
 
QED

\nocite{*} 
\bibliographystyle{unsrt}
\bibliography{bibli} %

\begin{thebibliography}{10}

\bibitem{ardentov2013conjugate}
AA~Ardentov and Yu~L Sachkov.
\newblock Conjugate points in nilpotent sub-riemannian problem on the engel
  group.
\newblock {\em Journal of Mathematical Sciences}, 195(3):369--390, 2013.

\bibitem{hakavuori2018blowups}
Eero Hakavuori and Enrico~Le Donne.
\newblock Blowups and blowdowns of geodesics in carnot groups.
\newblock {\em arXiv preprint arXiv:1806.09375}, 2018.

\bibitem{montgomery2002tour}
Richard Montgomery.
\newblock {\em A tour of subriemannian geometries, their geodesics and
  applications}.
\newblock Number~91. American Mathematical Soc., 2002.

\bibitem{bryant1986reduction}
Robert Bryant and Phillip Griffiths.
\newblock Reduction for constrained variational problems and∫ $\kappa$ 2/2
  ds.
\newblock {\em American Journal of Mathematics}, 108(3):525--570, 1986.

\bibitem{bonnard2001stratification}
Bernard Bonnard and Emmanuel Tr{\'e}lat.
\newblock Stratification du secteur anormal dans la sphere de martinet de petit
  rayon.
\newblock In {\em Nonlinear control in the Year 2000}, pages 239--251.
  Springer, 2001.

\bibitem{lawden2013elliptic}
Derek~F Lawden.
\newblock {\em Elliptic functions and applications}, volume~80.
\newblock Springer Science \& Business Media, 2013.

\bibitem{jurdjevic1997geometric}
Velimir Jurdjevic, Jurdjevic Velimir, and Velimir {\DJ}ur{\dj}evi{\'c}.
\newblock {\em Geometric control theory}.
\newblock Cambridge university press, 1997.

\bibitem{Bryant1993}
Robert Bryant and Hsu~Lucas L.
\newblock Rigidity of integral curves of rank 2 distributions.
\newblock {\em Inventiones mathematicae}, 114(2):435--462, 1993.

\bibitem{arnol2013mathematical}
Vladimir~Igorevich Arnol'd.
\newblock {\em Mathematical methods of classical mechanics}.
\newblock Springer Science, 1988.

\bibitem{marsden2013introduction}
Jerrold~E Marsden and Tudor~S Ratiu.
\newblock {\em Introduction to mechanics and symmetry: a basic exposition of
  classical mechanical systems}, volume~17.
\newblock Springer Science \& Business Media, 2013.

\bibitem{abraham1978foundations}
Ralph Abraham, Jerrold~E Marsden, and Jerrold~E Marsden.
\newblock {\em Foundations of mechanics}, volume~36.
\newblock Benjamin/Cummings Publishing Company Reading, Massachusetts, 1978.

\bibitem{bloch2003nonholonomic}
Anthony~M Bloch.
\newblock Nonholonomic mechanics.
\newblock In {\em Nonholonomic mechanics and control}. Springer, 2003.

\bibitem{hall2015lie}
Brian Hall.
\newblock {\em Lie groups, Lie algebras, and representations: an elementary
  introduction}, volume 222.
\newblock Springer, 2015.

\bibitem{montgomery2001geometric}
Richard Montgomery and Michail Zhitomirskii.
\newblock Geometric approach to goursat flags.
\newblock In {\em Annales de l'Institut Henri Poincare (C) Non Linear
  Analysis}, volume~18, pages 459--493. Elsevier, 2001.

\bibitem{stillwell2008naive}
John Stillwell.
\newblock {\em Naive lie theory}.
\newblock Springer Science \& Business Media, 2008.

\end{thebibliography}
\end{document}